\documentclass[11pt]{article}
\usepackage{graphicx}
\usepackage{amssymb}
\usepackage{amsmath}
\usepackage{amsthm}
\usepackage{color}
\numberwithin{equation}{section}
\newtheorem{theorem}{Theorem}

\newtheorem{cor}[theorem]{Corollary}

\newtheorem{oss}{Remark}

\newtheorem*{theorem*}{Theorem}

\title{Directional ellipticity on special domains: weak Maximum and Phragm\`en-Lindel\"of principles \vskip-0.25cm}
\author{Italo Capuzzo Dolcetta$^{\small 1}$, Antonio Vitolo$^{\small 2}$}
\date{$^{\small 1}$ \small Dipartimento di Matematica, Sapienza Universit\`{a} di Roma, 00185 Roma, ITALY
\vskip-0.425cm
$^{\small 2}$ Dipartimento di Ingegneria Civile, Universit\`{a} di Salerno, 84084 Fisciano (SA), ITALY, and\\ Istituto Nazionale di Alta Matematica, INdAM - GNAMPA
\vskip-0.375cm
}

\begin{document}
\maketitle

\noindent {\bf Abstract.} { We prove the validity of maximum principles for a class of fully nonlinear operators 
on unbounded subdomains $\Omega \subset \mathbb R^n$ of cylindrical type.
The main structural assumption is the uniform ellipticity of the operator
along the bounded directions of $\Omega$, with possible degeneracy along the unbounded directions.}



\section{Introduction}\label{INTRO}

The aim of this paper is to show that the weak Maximum Principle for $F$ in $\Omega$, that is 
\begin{equation}\nonumber
u \in USC(\overline \Omega),\ F(x,u,Du,D^2u)\geq 0 \  \hbox{\rm in} \ \Omega, \ u \ \le 0 \ \hbox{\rm on} \ \partial \Omega \quad \hbox{implies}\quad u \le 0  \  \hbox{\rm in} \ \Omega\, \eqno\hbox{\rm({\bf MP})}
\end{equation}
holds for a large class of degenerate elliptic nonlinear mappings $F$ and unbounded domains $\Omega$ of $\mathbb R^n$ whose geometry is related to the direction of ellipticity.\\
In the above, $\Omega$ is an open connected subset of $R^n$, $USC(\overline \Omega)$ denotes the set of upper semicontinuous functions on  $\overline \Omega$ and ${\mathcal S}^n$ is the set of the $n \times n$ real symmetric matrices with the usual partial ordering $X \le Y$ meaning that $Y-X$ is positive semidefinite.\\
For $u \in C^2(\Omega)$, $Du$ and $D^2u$ denoting, respectively, the gradient and the Hessian of the function $u$, the differential inequality in {\rm({\bf MP})} has the classical pointwise meaning.\\
On the other hand, for upper semicontinuous $u$ the partial differential inequality has to be understood in the viscosity sense, see  \cite{CIL},  \cite{CafCab}. 

\noindent We start describing the assumptions on the elliptic operator $F$, whose directions of strict ellipticity are related to the bounded directions of the domain $\Omega$  (pages 1 to 3), which will be all collected in condition  $(SC)_{U,\Omega}$. Next, we discuss a few examples of elliptic operators satisfying condition $(SC)_{U,\Omega}$, which are not necessarily uniformly elliptic in $\Omega$ (page 4 to 5).  Next, we present the main results (pages 5 to 8), and finally a comparison with the existing literature (pages 8 to 9).

\noindent Here $F=F(x,s,p,X)$  is a real-valued mapping defined in $\Omega\times\mathbb R\times\mathbb R^n\times {\mathcal S}^n$ 
which is assumed to satisfy
\begin{equation}\label{continuity} 
F \; \hbox{is continuous from}\; \Omega\times\mathbb R\times\mathbb R^n\times {\mathcal S}^n \;\hbox{into}\; \mathbb R\,,
\end{equation}
the monotonicity conditions
\begin{equation}\label{deg-elliptic} 
 F\; \mbox {is degenerate elliptic, that is}\;
 F(x,s,p,Y) \ge F(x,s,p,X) \quad \hbox{\rm if} \ Y\ge X
\end{equation}
and
\begin{equation}\label{s-monotonicity}
\hskip-2.75cm F(x,s,p,X)\le F(x,r,p,X)   \quad  \hbox{\rm if} \ s> r\,
\end{equation} 
as well as
\begin{equation}\label{F(0)=0}
\hskip-5.25cm F(x,0,0,O)=0 \quad \forall\,x \in \Omega\,
\end{equation} 
where $O$ is the zero-matrix.

In order to single out the class of domains $\Omega$ that we consider it is convenient to decompose $\mathbb R^n$ as the direct sum $\mathbb R^n=U\bigoplus U^\bot$, where  $U$ is a $k$-dimensional subspace and $U^\bot$ is its orthogonal complement. We shall denote for later purpose by $P$ and $Q$ the projection matrices on $U$ and $U^\bot$, respectively.\\
 We will assume then that  the open connected set $\Omega$  satisfies the following condition
\begin{equation}\label{slab-d}
\Omega \subseteq \{x \in \mathbb R^n : a \le x\cdot \nu^h \le a+d_h, \; h=1,\dots,k\}:=C\quad \mbox{for some} \; a \in \mathbb R, d_h>0\,,
\end{equation}
where $\{\nu^1,\dots,\nu^k\}$ is an orthonormal system for the subspace $U$.\\
Domains such as $C$, which we will sometimes refer to as $(n-k)$-infinite cylinders, with $1 \le k \le n-1$, for $n \ge2$,  are  infinite parallelepipeds whose $k$-dimensional orthogonal section is a $k$-parallelepiped of edges $d_h$, $h=1,\dots,k$. \\
In particular, a $(n-1)$-infinite cylinder is a slab, bounded in one direction and unbounded in all the remaining orthogonal directions.

\noindent It is worth pointing out here that such a domain $\Omega$  is typically unbounded, perhaps of infinite Lebesgue measure $| \Omega |$, but it does satisfy the measure-geometric condition {\bf(G)}, used in \cite{Cab} to obtain an improved form of the Alexandrov-Bakelman-Pucci ABP estimate, namely: \\
{\it there exist $R\in \mathbb R_+$, $\sigma,\tau \in (0,1)$ such that, for each $y \in \Omega$,  there is a ball $B_y=\{x\in \mathbb R^n : |x-y_0|<R_y\}$, with $R_y\le R$, providing the inequality $|B_y\backslash\Omega_{y,\tau}| \ge \sigma |B_y|$, where $\Omega_{y,\tau}$ is the connected component of $\Omega \cap B_{y,\tau}$ containing $y$, and $B_{y,\tau}=\{x\in \mathbb R^n : |x-y_0|<R_y/\tau\}$}.
\vskip-0.25cm
\noindent The above condition, first introduced in \cite{BNV}, requires, roughly speaking, that there is  enough boundary near every point in $\Omega$ allowing so to carry the information on the sign of $u$ from the boundary to the interior of the domain. It is therefore satisfied for example by unbounded domains of finite measure with $R= C(n) |\Omega|^{\frac{1}{n}}$  and also for a large class of unbounded domains with possibly infinite Lebesgue measure such as infinite cylinders, which we will be dealing with, or also perforated planes
$$ \mathbb R^2_{per}= \mathbb R^2\setminus \bigcup_{(i,j) \in \mathbb{Z}^2} B_{r}(i,j)$$
where $B_{r}(i,j)$ is the disc of radius $r<1$ centered at $(i,j)$. \\
Observe that $\sup_{y\in\mathbb R^2_{per}} \textrm{dist}(y,\partial{\mathbb R^2_{per}})<+\infty$, as for all domains satisfying condition {\bf(G)}. On the other hand, this fails to hold on cones, which therefore do not satisfy  {\bf(G)}\,.

\noindent It should also be stressed that no regularity assumption is made on the boundary $\partial \Omega$ so that the classical approach to establish comparison properties based on the construction of smooth barrier functions is not applicable in our general framework.

The next assumption is that there exists some $\nu \in\{\nu^1,\dots,\nu^k\}$ such that
\begin{equation}\label{elliptic-one-dir}
F(x,0,p,X+t\nu \otimes \nu)-F(x,0,p,X) \ge \lambda(x) t \quad \mbox{for all}\,\,t >0
\end{equation}
where  $\lambda$ is a continuous, strictly positive function such that $\liminf_{x\to \infty}\lambda(x)>0$.\\
Condition (\ref{elliptic-one-dir}) involves just the projection matrix $\nu \otimes \nu$ over the one dimensional subspace of $U$ spanned by the vector $\nu$. This strict ellipticity condition on $F$ related to the geometry of $\Omega$ will play a crucial role in our results.\\
We will assume moreover that there exist $\Lambda_1>0$  and a continuous function $\Lambda(x)\ge0$ such that $\Lambda(x)\le\Lambda_1|x|$ and
\begin{equation}\label{growth-unbd-dir}
F(x,0,0,X+tQ)-F(x,0,0,X) \le \Lambda(x)\,t \quad \mbox{for all}\,t >0, \mbox{as}\, |x|\to \infty
 \end{equation}
where $Q$ is the orthogonal projection matrix onto the subspace $U^\bot$. As for the behavior of $F$ with respect to the $p$ variable we assume that
\begin{equation}\label{growth-first-order}
 | F(x,0,p,X)-F(x,0,q,X)| \le \gamma(x) |p-q| \quad \mbox{for all}\,p,q\in \mathbb R^n
\end{equation}
with $\gamma(x)$ continuous, bounded  and such that $\frac{\gamma(x)}{\lambda(x)}$ is bounded above  in $\Omega$ by some constant $\Gamma\ge 0$.

\noindent We will refer collectively to conditions (\ref{continuity}), (\ref{deg-elliptic}), (\ref{s-monotonicity}), (\ref{F(0)=0}), (\ref{elliptic-one-dir}), (\ref{growth-unbd-dir}), (\ref{growth-first-order})
as the structure condition on $F$, labelled $(SC)_{U,\Omega}$. 

\noindent Observe that both $\nu \otimes \nu$ and $Q$ belong to ${\mathcal S}^n$ and are positive semidefinite. It is also worth noting that conditions (\ref{elliptic-one-dir}) and (\ref{growth-unbd-dir}), requiring, respectively, a control from below only with respect to a single direction $\nu\in U$ and a control from above in the orthogonal directions, comprise a much weaker condition on $F$ than uniform ellipticity.\\
The latter one would indeed require a uniform control of the difference quotients both from below and from above with respect to all possible increments with positive semidefinite matrices.

\noindent Consider for simplicity the case  $U=\{x=(x_1,...,x_n)\in \mathbb R^n: x_{k+1}=\dots=x_n=0\}$ and $\Omega=\{x \in \mathbb R^n : 0 \le x_h \le d, \; h=1,\dots,k\}$, where $d>0$.   A very basic example of  an $F$ satisfying $(SC)_{U,\Omega}$ is given by the linear operator
\begin{equation}\label{ex1:linear}
F(x,u,Du,D^2u)=\sum_{i=1}^k\frac{\partial^2u}{\partial x_i^2}+|x| \sum_{i=k+1}^{n}\frac{\partial^2u}{\partial x_i^2}+ \sum_{j=1}^nb_j(x)\, \frac{\partial u}{\partial x_j} + c(x)u
\end{equation}
with $b_i(x)$ and $c(x)$ continuous functions such that 
\begin{equation}\label{lot}
\left|\sum_ib_i^2(x)\right|^{1/2}\le \gamma < \infty, \;\;\; c(x) \le 0.
\end{equation}
In fact, $F(x,s,p,X)= \text{Tr}(AX) + b\cdot p+cs$, where
$$A= \left(\begin{array}{cc}
\mathbb I_k & \mathbf 0 \hskip0.3cm\\
\mathbf 0 &  |x| \, \mathbb I_{n-k}
\end{array}\right), \;\;\; b=(b_1,\dots,b_n),$$
which of course satisfies (\ref{continuity}), (\ref {deg-elliptic}) and (\ref{F(0)=0}). \\
Next, (\ref {s-monotonicity}) is ensured by the sign condition on $c(x)$ in (1.10). Moreover, looking at the above matrix $A(x)$, it is immediate to check that $F$ is uniformly elliptic with respect to any direction of $U$, so verifying (\ref{elliptic-one-dir}). On the other hand, the differential quotients with respect to matrix increments related to the orhogonal directions $x_{k+1}, \dots, x_n$ are $O(|x|)$, which imply (\ref{growth-unbd-dir}). The assumption (\ref{lot}) on $b_i(x)$ provides  (\ref{growth-first-order}).
Hence we conclude that actually $F$ in (\ref{ex1:linear}) satisfies $(SC)_{U,\Omega}$.

\noindent Relevant nonlinear examples are provided by fully nonlinear operators of Bellman-Isaacs type
\begin{equation}\label{ex1:fullynonlinear}
F(x,u,Du,D^2u)=\sup_{\alpha}\inf_\beta L^{\alpha\beta}u,
\end{equation}
where 
$$L^{\alpha\beta}u = \sum_{i,j=1}^na_{ij}^{\alpha\beta}\frac{\partial^2u}{\partial x_ix_j} + \sum_{i=1}^nb_i^{\alpha\beta}\frac{\partial u}{\partial x_i} + c^{\alpha\beta} u$$ 
with constant coefficients depending $\alpha$ and $\beta$ running in some sets of indexes $\cal A, \cal B$. \\
If $A^{\alpha\beta}=[a_{ij}^{\alpha\beta}]$ is  positive semidefinite for all $\alpha,\beta$ and 
\begin{align*}
&\sum_{i,j=1}^n a_{ij}^{\alpha\beta}\nu^h_i\nu^h_j\ge \lambda,  \;\;\;  h=1,\dots,k,  \\
&\sum_{i,j=1}^n a_{ij}^{\alpha\beta}\nu^h_i\nu^h_j a_{ij}^{\alpha\beta}\nu^h_i\nu^h_j \le \Lambda|x|,  \;\;\;  h=k+1,\dots,n \\ 
&|b_i^{\alpha\beta}| \le \gamma, \quad c^{\alpha\beta}\le 0,
\end{align*}
where  $\{\nu^1,\dots,\nu^k\}$ is an orthonormal basis of $U$, then $F$ satisfies our assumptions in any domain $\Omega$ contained in a $(n-k)$-infinite cylinder like (\ref{slab-d}).

\noindent Our results concerning the validity of ({\bf MP}) are stated in the following theorems:
\begin{theorem}\label{MP-slab-thm}  Let $\Omega$ be a domain of $\mathbb R^n$ satisfying condition $(\ref{slab-d})$ and assume that $F$ satisfies $(SC)_{U,\Omega}$.
Then ({\bf MP}) holds for any $u \in USC(\overline \Omega)$ such that $u^+_0=o(|x|)$ as $|x| \to \infty$.
\end{theorem}
\noindent The function  $u^+_0$ in the statement is defined by
\begin{equation*}
u^+_0(x)=\left\{\begin{array}{ll}
u^+(x)\equiv\max(u(x),0) & \hbox{\rm if} \ x \in \overline\Omega\\
0 & \hbox{\rm if} \ x \not\in \overline\Omega
\end{array}\right.
\end{equation*}
and condition  $u^+_0(x)=o(|x|)$ as $|x| \to \infty$ is understood as 
$$\limsup_{\substack{x\in\Omega\\ x \to \infty}}\frac{u^+_0(x)}{|x|}=0.$$
 The use of $u^+_0$ is convenient for a unified statement of our results which are valid both for bounded and unbounded $\Omega$, with the obvious exception of  Theorem \ref {Phragmen} below on Phragm\`en-Lindel\"of principles.

\noindent Note that some restriction on the behaviour of $u$ at infinity is unavoidable. Observe in this respect that, in the 1-infinite cylinder $C_1=\mathbb R\times (0,\pi)^2\subset \mathbb R^3$, the function $u(x_1,x_2,x_3)=e^{x_1}\sin x_2\sin x_3$  solves the degenerate Dirichlet problem
 $$
\frac{\partial^2 u}{\partial x_1^2} +\frac{\partial^2 u}{\partial x_2^2} = 0 \ \mbox{in}\ \ C_1, \ \  u(x_1,x_2,x_3) = 0\,\;  \mbox{on}\,\; \partial C_1
$$ 
 and $u(x_1,x_2,x_3)>0$ in $C_1$, implying the failure of ({\bf MP}).\\
The above condition  $u^+(x)=o(|x|)$ at infinity is in accordance, by the way, with the assumption of Ishii in \cite[Section 7]{I}, where a comparison principle in unbounded domains is obtained under the assumption of at most linear growth of the solutions at infinity.

\noindent Let us consider now the degenerate elliptic operator 
$$Lu=\lambda_1(x)\,\frac{\partial^2 u}{\partial x_1^2} +\lambda_2(x)\,\frac{\partial^2 u}{\partial x_2^2}+\lambda_3(x)\,\frac{\partial^2 u}{\partial x_3^2}$$ 
where $\lambda_i(x)$ are continuous functions such that $\lambda_i(x) >0$, $i=1,2,3$, such $\lambda_i(x)\le\Lambda_1|x|$, in the $1$-infinite cylinders 
$$C_1=\mathbb R\times (0,\pi)^2, \;\;\;C_2=(0,\pi)\times \mathbb R\times (0,\pi), \;\;\;C_3=(0,\pi)^2\times \mathbb R,$$ 
which have $x_i$ as unbounded direction. 
Therefore, according  to Theorem \ref{MP-slab-thm}, provided 
$$\inf_{C_i} \lambda_j(x)>0 \;\;\; \mbox{for}\; j\neq i$$ 
and $u^+(x)=o(|x|)$, this operator satisfies  ({\bf MP}) in all domains contained in $C_1$ as well in $C_2$ and in $C_3$. As a further application,  noting that the intersection $C_1\cap C_2\cap C_3$ is the cube $(0,\pi)^3$, where $L$ is uniformly elliptic, we infer that ({\bf MP}) holds also in all domains contained in $C_1\cup C_2\cup C_3$.

\noindent More generally, below we are going to show a consequence of Theorem \ref{MP-slab-thm}, which supports this claim as a particular case.\\
As above, let  $\{\nu^1,\dots,\nu^n\}$ an orthonormal basis of $\mathbb R^n$, and let $C^{k}$  be  $(n-1)$-cylinders having one of $\nu^i$, $i=1,\dots,n$ as axis. We consider lattice domains which are finite union of  $1$-infinite cylinders. See Fig.\ref{lattice} below.

\begin{figure*}[ht]
\hskip2.0cm
\includegraphics[height=6cm,width=12cm,angle=0]{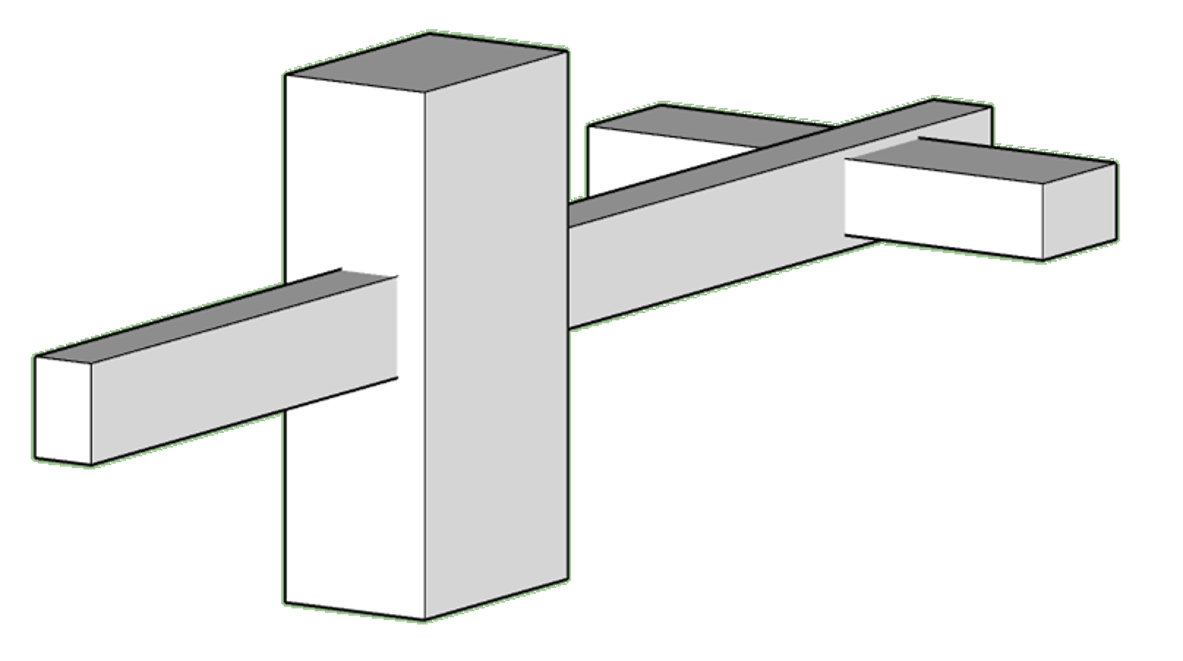}
\caption{Spatial lattice}
\label{lattice}
\end{figure*}
\noindent Our approach to prove this claim makes use of the strong Maximum Principle, according to which a subsolution of a uniformly elliptic equation cannot have a non-negative maximum inside the domain unless it is constant. See \cite{GT} for linear elliptic operators and  {\cite{BDL}, \cite{ARV} for different approaches in the nonlinear case.\\
This explains  why in the next result the uniform ellipticity will be assumed in the nodes where two or more cylinders intersect each other, which is not guaranteed by the one-directional ellipticity with respect to orthogonal directions occurring in the nodes as regions belonging to $1$-infinite cylinders with different axes.
\begin{cor}\label{cross} Suppose  $\Omega$ to be a domain contained in a lattice domain $R$ of $1$-infinite cylinders $C$ with axes $\nu$ in an orthonormal basis $\{\nu^1,\dots,\nu^n\}$ of $\mathbb R^n$. Letting $U$ be the $(n-1)$- dimensional  subspace orthogonal to the axis $\nu$, we also suppose that $F$ satisfies the structure condition $(SC)_{U,C}$ in each $C \subset R$.\\
Let also $N$ be the union of points belonging to more than one cylinder. If in addition $F$ is uniformly elliptic in $N$,  then $ ({\bf MP})$ holds in $\Omega$, provided $u^+_0(x)=o(|x|)$ as $|x|\to\infty$.
\end{cor}
\noindent The next result can be seen as a quantitative form of Theorem \ref {MP-slab-thm} above, it is an $L^\infty$ version of the ABP estimate.
\begin{theorem}\label{ABP-simple-thm} Let $\Omega$ be a domain of $\mathbb R^n$ satisfying condition $(\ref{slab-d})$ and assume that $F$ satisfies $(SC)_{U,\Omega}$.
If $u \in USC(\overline \Omega)$ satisfies $u^+_0(x)=o(|x|)$ as $|x| \to \infty$ and
$$F(x,u,Du,D^2u) \ge f(x)\quad \hbox{in}\; \Omega\,$$ 
where $f/\lambda$ is continuous and bounded from below, then  for all $h=1,\dots,k$
\begin{equation}\label{simple-ABP-est}
\sup_{\overline\Omega} u \le \sup_{\partial\Omega}u^+ + \frac{e^{1+d_h\,\Gamma }}{1+d_h\,\Gamma} \left\|\frac{f^-}{\lambda}\right\|_\infty d_h^2
\end{equation}
where $f^-(x)=-\min(f(x),0)$. 
\end{theorem} 

\noindent The next result is about bounded functions $u$ in the so-called narrow domains, that is when $(\ref{slab-d})$ is fulfilled with at least one $d_h$ sufficiently small. In this case we may assume that $F$ is just nondecreasing with respect to the variable $s$:
\begin{theorem}\label{narrow-thm}  Let $\Omega$ be a domain of $\mathbb R^n$ satisfying condition $(\ref{slab-d})$ and assume that $F$ satisfies $(SC)_{U,\Omega}$ with  (\ref{s-monotonicity}) replaced by the weaker condition
\begin{equation}\label{s-bound-derivative}
 F(x,s,p,M)-F(x,r,p,M) \le c(x)\,(s-r)    \ \ \hbox{\rm if} \ s>r
\end{equation}
  for some continuous function $c(x)>0$. Assume also that  $\frac{c(x)}{\lambda(x)}\leq K<+\infty$  in $\Omega$ with $\lambda$ as in  (\ref {elliptic-one-dir}). Then $({\bf MP})$ holds for $u \in USC(\overline \Omega)$,  $u$ bounded above, provided $ d^2 K$ is small enough, where $d=\min_{h}d_h$ in $(\ref{slab-d})$.
\end{theorem}

\noindent The above result is often used as an intermediate step in the proof of Theorem \ref{Phragmen} below concerning the validity of $({\bf MP})$ for unbounded solutions with exponential growth at infinity.
\begin{theorem}\label{Phragmen}
 Let $\Omega$ be a domain of $\mathbb R^n$ satisfying condition $(\ref{slab-d})$ and assume that $F$ satisfies $(SC)_{U,\Omega}$. In addition, assume that there exists a positive number $\rho$ such that $\frac{\Lambda(x)}{\lambda(x)}\le \rho <\infty$. \\
 Then, for any fixed $\beta_0 >0$ there exists a positive constant $d=d(n,\lambda,\Lambda,\gamma,\beta_0)$ such that if $d_h \le d$ for some $h=1,\dots,k$, then $({\bf MP})$ holds for functions $u$ such that $u^+_0(x)=O(e^{\beta_0|x|})$ as $|x| \to \infty$.\\
Conversely, for any fixed $d_0 >0$, supposing $d_h \le d_0$ for all $h=1,\dots,k$, there exists  a positive constant $\beta=\beta(n,\lambda,\Lambda,\gamma,d_0)$ such that $({\bf MP})$ holds for functions $u$ such that $u^+_0(x)=O(e^{\beta|x|})$ as $|x| \to \infty$.
\end{theorem}

\noindent Note that the growth control (\ref{growth-unbd-dir})  in the directions of $U^\bot$ is essential in order to have the maximum principle for subsolutions growing at infinity more than polynomially, as the following example shows. \\
Indeed, $u(x_1,x_2)=x_2^2\sin x_1$ is a solution of $$
\frac{\partial^2u}{\partial x_1^2}+\frac12\,x_2^2\,\frac{\partial^2u}{\partial x_2^2}=0
$$
in the cylinder $C= (0,\pi) \times \mathbb R \subset \mathbb R^2$, $u=0$ on $\partial\Omega$ but $u$ is strictly positive in $C$ .

\noindent Maximum and Phragm\`en-Lindel\"of principles  are extensively dealt with in the book of Protter and Weinberger \cite{PW} for classical solutions of linear uniformly elliptic operator. A first Maximum Principle for strong solutions in nonsmooth domains of cylindrical type was proved by Cabr\'e \cite{Cab} using  a measure-geometric condition, which originates from a work of Berestycki, Nirenberg and Varadhan \cite{BNV}. Such a condition was generalized in \cite{CV} and \cite{V} to include a larger class of domains like conically and parabolically shaped domains, and extended to the fully nonlinear setting in the viscosity sense in \cite{CDLV}. Results on the weak MP in cylindrical domains have been also established by Busca \cite{Bu}.\\
Further results on ({\bf MP}) in unbounded domains with superlinear gradient terms as well as Phragm\`en-Lindel\"of principles  have been obtained in \cite{ARV} and in \cite{CDV} for the natural quadratic growth in the gradient. For further Phragm\`en-Lindel\"of principles see also \cite{Pun}. \\
More recently, maximum principles in domains of cylindrical type as well as global H\"older estimates have been shown for degenerate or singular elliptic operators of $p$-Laplacian type have been proved in \cite{BD}, \cite{BCDV}, based on the results of \cite{Im}. Different maximum principles for subsolutions  which may be unbounded at finite points and their application to removable singularities issues for degenerate elliptic operators of different type, like partial sums of eigenvalues, have been stated in \cite{GV} and \cite{VIT}. \\
The most closely related paper to our present work is \cite{CDV2}, where maximum and Phragm\`en-Lindel\"of principles  for viscosity subsolutions in domains of cylindrical type $\Omega$ which are bounded just in one single direction, called $(n-1)$-infinite cylinders, assuming a strict ellipticity in the bounded direction and a bound from above on the difference quotients of $F$ with respect to the increments of the matrix variables in the orthogonal directions. Existence results in such domains can be found in \cite{Man}.\\
In the present paper, a weaker control is assumed on the differential quotients of $F$, which have no more to be bounded in the unbounded directions (at most linear growth with respect to $x$ is admissible). Moreover, the case of more bounded directions is considered, covering the case of $(n-k)$-infinite cylinders with $1<  k<n$. The advantage of this sort of multidimensional boundedness lies in the possibility of assuming the control on the differential quotients with respect to less directions, more precisely $(n-k)$ with $k>1$ instead of $(n-1)$. On the other hand, the role of strict ellipticity in more than one direction is pointed out in particular by considering domains which are union of $1$-infinite cylinders.

\section{The weak Maximum Principle via one-directional strict ellipticity}\label{MP-strict-ellipticity} 
This section is devoted to the proof of Theorem \ref{MP-slab-thm}  concerning the validity of the weak Maximum Principle in domains contained in $(n-k)$-infinite cylindrical domains.\\
The next simple observation exploits a useful consequence  of condition (\ref{F(0)=0}) and the growth condition (\ref{growth-unbd-dir}) with respect to the so-to-say unbounded directions:\\
\begin{equation}\label{technical*}\limsup_{\varepsilon \to 0^+} F(x_\varepsilon,0,0,\tfrac{\varepsilon}{|x_\varepsilon|}\,Q) \le 0\; \mbox {for any sequence}\,  x_\varepsilon\in\Omega\, \mbox{such that} \lim _{\varepsilon \to 0^+}|x_\varepsilon| = +\infty 
\end{equation}

\noindent Indeed, if $x_\varepsilon\in \Omega$ is such a sequence, an immediate consequence (\ref{growth-unbd-dir}) with $X=O$ and $t=\tfrac\varepsilon{|x_\varepsilon|}$ is that 
$$F(x_\varepsilon,0,0,\tfrac\varepsilon{|x_\varepsilon|}\,Q)= F(x_\varepsilon,0,0,\tfrac\varepsilon{|x_\varepsilon|}\,Q)- F(x_\varepsilon,0,0,O)\\
\le \Lambda_1 \, \varepsilon
$$

{\bf Proof of Theorem \ref {MP-slab-thm} } 
We may assume that $U=\mathbb R^k \times \{0\}^{n-k}$ so that $U^\bot=\{0\}^{k}\times \mathbb R^{n-k}$ and
$$
P=\left(\begin{array}{cl}
\mathbb I_k & 0\\
0 & 0_{n-k}
\end{array}\right), \quad \quad 
Q= \left(\begin{array}{cl}
0_k & 0\\
0 & \mathbb I_{n-k}
\end{array}\right),
$$
$0_{k}$ and $\mathbb I_{k}$ being the $k \times k$ zero and identity matrices, respectively.\\
We may also assume that  $\sqrt{x_1^2+\dots+x_k^2} \le d$ for some $d\in\mathbb R_+$ for all $x=(x_1,\dots,x_n) \in \Omega$.\\
Arguing by contradiction, suppose that $u(x)$ has a positive value $M>0$ at some point $x \in \Omega$. 
For $\varepsilon >0$ we consider the function 
\begin{equation}\label{ueps}
u_\varepsilon(x)=u(x)-\varepsilon \varphi(x) 
\end{equation}
where $\varphi(x)=\sqrt{x_{k+1}^2+\dots+x_n^2+1}$. 
Since $u(x)=o(|x|)$ as $|x|\to \infty$, then $u_\varepsilon(x)\le0$ for $|x|$ large enough. So there exists a bounded domain $\Omega_\varepsilon \subset \Omega$ such that $u_\varepsilon(x)\le 0$ for $x \not\in\overline\Omega_\varepsilon$, and
 \begin{equation}\label{Meps}
M_\varepsilon\equiv\sup_{\Omega}u_\varepsilon =  \max_{\overline\Omega_\varepsilon}u_\varepsilon =u(x_\varepsilon)\ge \frac M2.
\end{equation}
Since $u_\varepsilon \le 0$ on $\partial\Omega_\varepsilon$ then $x_\varepsilon \in \Omega_\varepsilon$.
On the other hand, since $u=u_\varepsilon+\varepsilon \varphi(x)$ is a subsolution, then
$u_\varepsilon$ satisfies the differential inequality
\begin{equation}\label{ueps:eq}
F(x,u_\varepsilon+\varepsilon \varphi(x),Du_\varepsilon+\varepsilon D\varphi(x),D^2u_\varepsilon+\varepsilon D^2\varphi(x)) \ge 0 \,.
\end{equation}
Let $\nu^h\equiv \nu=\nu'\times 0_{n-k}$ be the direction of uniform ellipticity in $U$, with $\nu'=(\nu_1',\dots,\nu_k')$, and $d$ the corresponding $d_h$ in (\ref{slab-d}). \\
 We consider, for $\alpha >0$ to be chosen in the sequel, the function
\begin{equation}
h_\varepsilon(x)=M_\varepsilon(1+e^{-\alpha d}) -M_\varepsilon e^{\alpha (x'-x_\varepsilon')-\alpha d},
\end{equation}
where $x'=\nu_1x_1+\dots+\nu_kx_k$ and $x_\varepsilon$ is the maximum point of (\ref{Meps}).\\
Note that  $h_\varepsilon(x_\varepsilon)=M_\varepsilon=u_\varepsilon(x_\varepsilon)$ and 
\begin{equation}\label{heps:ineq}
h_\varepsilon \ge M_\varepsilon e^{-2\alpha d} \ge \frac M2\,e^{-2\alpha d}\quad \mbox{ in}\  \overline\Omega.
\end{equation}
Modulo the addition of some number $c_\varepsilon \ge 0$,  we can make $h_\varepsilon+c_\varepsilon$ touch $u_\varepsilon$ from above at a point $\overline x_\varepsilon \in \Omega_\varepsilon$. Therefore  $h_\varepsilon(x)+c_\varepsilon+\varepsilon\varphi$  can be used as a test function for the subsolution $u$, and will  satisfy the inequality 
\begin{equation}\label{F(heps)}
F(\overline x_\varepsilon,h_\varepsilon(\overline x_\varepsilon)+c_\varepsilon+\varepsilon \varphi(\overline x_\varepsilon),Dh_\varepsilon(\overline x_\varepsilon)+\varepsilon D\varphi(\overline x_\varepsilon),D^2h_\varepsilon(\overline x_\varepsilon)+\varepsilon D^2\varphi(\overline x_\varepsilon)) \ge 0.
\end{equation}
Since $h_\varepsilon(\overline x_\varepsilon)+c_\varepsilon+\varepsilon \varphi(\overline x_\varepsilon)>0$ then by monotonicity condition (\ref{s-monotonicity})
\begin{equation}\label{eq1-MP-slab-thm}
F(\overline x_\varepsilon,0,Dh_\varepsilon(\overline x_\varepsilon)+\varepsilon D\varphi(\overline x_\varepsilon),D^2h_\varepsilon(\overline x_\varepsilon)+\varepsilon D^2\varphi(\overline x_\varepsilon))\ge 0\,.
\end{equation}

\noindent Computing the derivatives  
\begin{align*}
&Dh_\varepsilon(\overline x_\varepsilon)= -\alpha M_\varepsilon e^{\alpha (\overline x_\varepsilon'-x_\varepsilon')-\alpha d} \,\nu, \\
&D^2h_\varepsilon(\overline x_\varepsilon)= -\alpha^2M_\varepsilon e^{\alpha (\overline x'_\varepsilon-x_\varepsilon')-\alpha d}\,\nu\otimes \nu,
\end{align*}
and observing that
$$M_\varepsilon e^{\alpha (x_1-x_{\varepsilon,1})-\alpha d} \ge  \frac M2\,e^{-2\alpha d},
$$
from  (\ref{elliptic-one-dir}) and (\ref{growth-first-order}) we get
\begin{equation}\label{eq2-MP-slab-thm}
\begin{split}
&F(\overline x_\varepsilon,0,Dh_\varepsilon(\overline x_\varepsilon)+\varepsilon D\varphi(\overline x_\varepsilon),D^2h_\varepsilon(\overline x_\varepsilon)+\varepsilon D^2\varphi(\overline x_\varepsilon))\\
\le & F(\overline x_\varepsilon, 0,\varepsilon D\varphi(\overline x_\varepsilon),{\varepsilon} D^2\varphi(\overline x_\varepsilon)) -\alpha \, M_\varepsilon e^{\alpha (x'-x'_\varepsilon)-\alpha d}(\lambda(\overline x_\varepsilon) \alpha - \gamma(\overline x_\varepsilon) )\\
\le& F(\overline x_\varepsilon, 0,\varepsilon D\varphi(\overline x_\varepsilon),{\varepsilon} D^2\varphi(\overline x_\varepsilon)) - C\alpha\lambda(\overline x_\varepsilon)\left(\alpha - \frac{\gamma(\overline x_\varepsilon)}{\lambda(\overline x_\varepsilon)} \right) 
\end{split}
\end{equation}
with $C= \frac12\,Me^{-2\alpha d}>0$.

From  (\ref{eq1-MP-slab-thm}) and (\ref{eq2-MP-slab-thm}), taking into account that $\frac{\gamma}{\lambda}\le \Gamma$, it follows that
\begin{equation}\label{eq3-MP-slab-thm}
\begin{split}
0 &\le F(\overline x_\varepsilon, 0,\varepsilon D\varphi(\overline x_\varepsilon),{\varepsilon} D^2\varphi(\overline x_\varepsilon) ) -C\alpha\lambda(\overline x_\varepsilon)\left(\alpha - \Gamma \right)\\
& \le  F(\overline x_\varepsilon, 0,\varepsilon D\varphi(\overline x_\varepsilon),{\varepsilon} D^2\varphi(\overline x_\varepsilon) ) -\tfrac12\,C\alpha^2\lambda(\overline x_\varepsilon)
\end{split}
\end{equation}
choosing $\alpha >2\Gamma$.

Now we consider the two possible cases:

(i) $\{\overline x_\varepsilon\}$ bounded; (ii) $\{\overline x_\varepsilon\}$ unbounded.

In case (i),  we can extract a subsequence converging to $\overline x \in  \Omega$. In fact, we may esclude that  $\overline x \in  \partial \Omega$, where $u \le 0$,  since by construction, using the upper semicontinuity of $u$, we have
$$u(\overline x) \ge \lim_{\varepsilon\to 0}u(\overline x_\varepsilon)\ge \frac M2>0.$$
Recalling that $F(\overline x,0,0,0)=0$, by the continuity of $F$ we have
$$
\lim_{\varepsilon \to 0}F(\overline x_\varepsilon, 0,\varepsilon D\varphi(\overline x_\varepsilon),{\varepsilon} D^2\varphi(\overline x_\varepsilon)) = 0,
$$
so that, also using the continuity of $\lambda$, from  (\ref{eq3-MP-slab-thm}) as $\varepsilon \to 0$  we get the following contradiction:
\begin{equation}\label{eq4-MP-slab-thm}
\tfrac12\,C\alpha^2\lambda(\overline x)\le 0,
\end{equation}
whereas $\lambda(\overline x)>0$  by assumption, and this concludes the proof in case (i).

In case (ii), where $\overline x_\varepsilon$ is unbounded, we take a subsequence such that $|\overline x_\varepsilon| \to \infty$ as $\varepsilon \to 0$. Computing the derivatives of $\varphi(x)$, we get
\begin{equation}\label{Dephi}\begin{split}
|D\varphi(x)| \le 1, \quad  D^2\varphi(x) \le  \frac{Q}{\varphi(x)}\,.
\end{split}
\end{equation}
From this, using (\ref{growth-first-order}) and degenerate ellipticity, we get
\begin{equation}\label{eq5-MP-slab-thm}
\begin{split}
F(\overline x_\varepsilon, 0,\varepsilon D\varphi(\overline x_\varepsilon),{\varepsilon} D^2\varphi(\overline x_\varepsilon)) 
&\le  F(\overline x_\varepsilon, 0,0,{\varepsilon}\, \tfrac{Q}{\varphi(\overline x_\varepsilon)}) +\varepsilon\,\sup_\Omega\gamma(x)
\end{split}
\end{equation}
so that, using  (\ref{technical*}), 
$$
\limsup_{\varepsilon \to 0}F(\overline x_\varepsilon, 0,\varepsilon D\varphi(\overline x_\varepsilon),{\varepsilon} D^2\varphi(\overline x_\varepsilon)) \le 0.
$$
Finally,  letting $\varepsilon \to 0^+$,  we estimate (\ref{eq3-MP-slab-thm}) with (\ref{eq5-MP-slab-thm}), so we again obtain a contradiction:
$$
\liminf_{\varepsilon\to 0}\lambda(\overline x_\varepsilon)\le 0,
$$
whereas $\liminf_{x\to \infty}\lambda(x)>0$  by assumption, concluding the proof.\qed

\noindent Here below we see how Corollary \ref{cross} follows from Theorem \ref{MP-slab-thm}.

\noindent{\bf Proof of Corollary \ref{cross}.} For sake of clarity, and without loss of generality, we illustrate the proof in the case $n=2$ with the union $R$ of two $1$-infinite cylinders of axes $(1,0)$ and $(0,1)$. Suppose $F[u]\ge 0$ in $\Omega$ and $u\le 0$ on $\partial \Omega$.  We may also suppose, eventually passing to $u^+_0$, that $\Omega=R$, $u\ge0$ in $R$ and $u=0$ on $\partial R$. To be more direct, we refer to Figure \ref{R2} below.
\begin{figure*}[ht]
\hskip2.0cm
\includegraphics[height=10cm,width=12cm,angle=0]{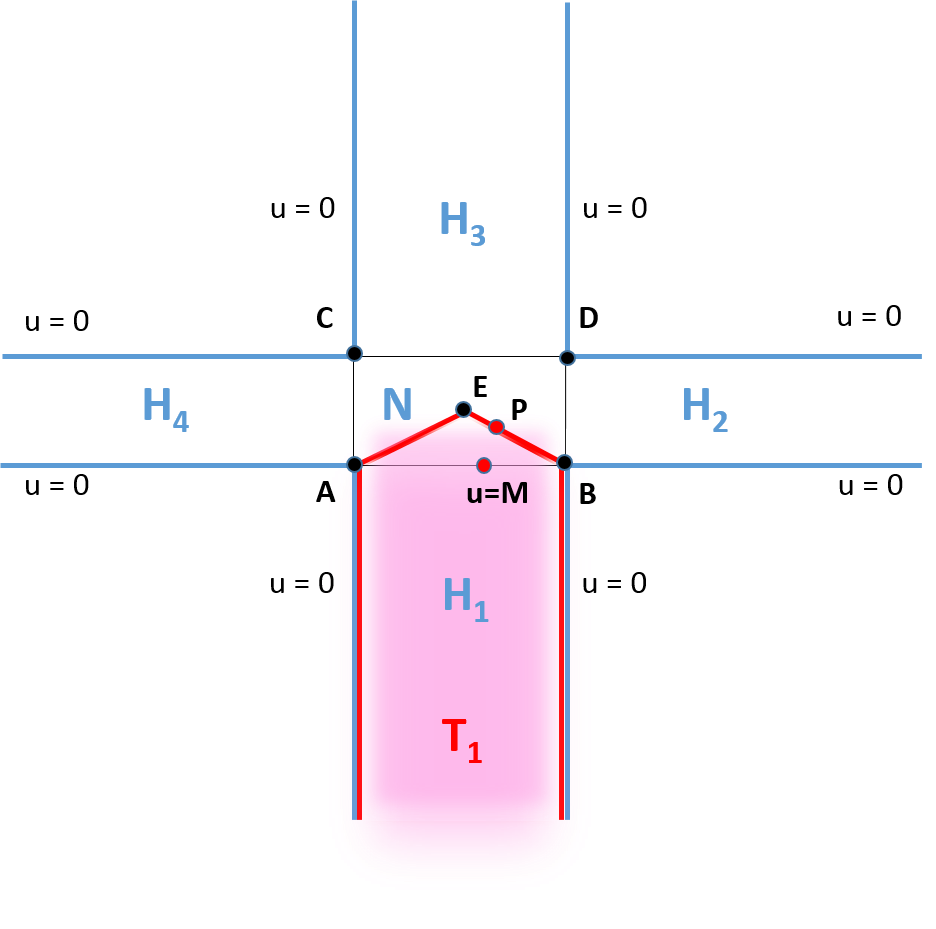}
\vskip-0.5cm
\caption{Crossing strips}
\label{R2}
\end{figure*}

\noindent  By the weak Maximum Principle of Theorem \ref{MP-slab-thm}, the maximum of $u$ on each half-strip $H_i$, $i=1,\dots,4$, is achieved on the boundary of the central rectangle $N$ as well as the maximum on $N$ itself, so that the maximum $M$ on the whole $R$ is on $\partial N$.\\
Again looking at Figure \ref{R2}, we suppose that $M$  is attained on the side $AB$. Considering the cylindrical domain $T_1$ (shaded area)  obtained from $S_1$ with the addition of the triangle $ABE$, then again by the weak Maximum Principle there is a maximum point $P$ on at least one of sides $AE$ and $BE$. But $P \in N$, where $F$ is assumed to be uniformly elliptic. Therefore, by the strong Maximum Principle, $u$ has to be constant in $R$, and this constant has to be necessarily zero, as it was to be shown.\hfill$\Box$

 \medskip
 A comparison principle between upper semicontinuous subsolutions $u$ and smooth supersolutions $v$
follows immediately  from Theorem \ref{MP-slab-thm} and Corollary \ref{cross} above. Here is the result:
\begin{cor}
Assume on $F$ and $\Omega$  the same conditions as in Theorem \ref {MP-slab-thm} or  Corollary \ref{cross}.
If $u \in USC(\overline \Omega)$ and $v \in C^2(\Omega)\cap LSC(\overline \Omega)$ satisfy
$$
F(x,v,Dv,D^2v) \le F(x,u,Du,D^2u) \,\, \mbox{\rm in}\,\ \Omega
$$
and  $(u-v)^+(x)=o(|x|)$ as $|x| \to \infty$, 
then
$$
u \le v \, \hbox{\rm on}\, \  \partial \Omega \quad \hbox{\rm implies} \quad u \le v \,\, \hbox{\rm in}\, \ \Omega\,\eqno
$$
\end{cor}
\noindent To prove this, it is enough to observe that the operator $G$ defined as 
$$G(x,s,p,M)=F(x,s+v(x),p+Dv(x),M+D^2v(x))-F(x,v(x),Dv(x),D^2v(x))$$
fulfills all the assumptions of Theorem \ref {MP-slab-thm}, noting that the statement is equivalent to the maximum principle for the subsolution $u-v$ of equation $G(x,w,Dw,D^2w)=0$.\quad $\Box$

\section {The weak Maximum Principle in narrow domains} 
To deal with the ({\bf MP}) in narrow domains we need the $L^\infty$ estimate of Theorem \ref{ABP-simple-thm}, a sort of quantitative version of ({\bf MP}), which can be deduced from  Theorem 1.
\noindent Here below we show a proof of this result.

\noindent {\bf Proof of Theorem \ref{ABP-simple-thm}.}  Let $u \in USC(\overline \Omega)$  be a viscosity solution of the inequality $F(x,u,Du,D^2u) \ge f(x)$ in $\Omega$. Conditions (\ref{F(0)=0}), (\ref{s-monotonicity}) and some viscosity calculus show that $u^+={\rm max}(u,0)$ satisfies
 \begin{equation}\label{subsol.f}
F(x,0,Du^+,D^2u^+) \ge -f^-(x)\,.
\end{equation}

We may assume that the direction $\nu$ in the ellipticity condition (\ref{elliptic-one-dir}) is the positive direction of the $x_1$ axis and that $0\le x_1 \le d_h$ for $x \in \Omega$. \\
Suppose temporarily that $d_h=1$ and consider the auxiliary function
$$
w(x)= u^+(x)+C_1e^{\alpha x_1} - \sup_{\partial\Omega}u^+ - C_1e^\alpha\,
$$
where $C_1$ and $\alpha$ are positive constants to be chosen in the sequel. \\
The structure condition (SC)$_{U,\Omega}$ yields
\begin{align*}
F(x,0,Dw,D^2w) &\ge F(x,0,Du^+,D^2u^+) -\gamma(x) C_1 \alpha e^{\alpha x_1} + \lambda(x) C_1 \alpha^2 e^{\alpha x_1}\\
&\ge  -f^-(x) + \alpha  C_1  (\alpha\lambda(x) -\gamma(x)) 
 \, \hbox{ in} \,\ \Omega\,
\end{align*}
Choosing $\alpha = 1+\Gamma$ with $\Gamma$ as in (\ref{growth-first-order}) and
$C_1=\frac{1}{1+\Gamma}\,\sup_{\Omega}\frac{f^-}{\lambda}$
we obtain
$$F(x,0,Dw,D^2w) \ge 0 \ \, \hbox{in}\,\ \Omega, \ \ w \le 0 \, \hbox{on}\,\partial\Omega\,.$$
By Theorem 1 we conclude that $w \le 0$ in $\Omega$, which implies
$$
u(x) \le u^+(x) + C_1e^{\alpha x_1}\le \sup_{\partial\Omega}u^+ + \frac{e^{1+\Gamma }}{1+\Gamma } \,\sup_{\Omega}\frac{f^-}{\lambda}\,\,
$$
proving (\ref{simple-ABP-est}) in the case $d=1$.\\
For an arbitrary $d_h >0$, we consider the rescaled variable $y= x/d_h$ and the operator $G$ defined by
$$
G(y,s,p,M) = d_h^{2}\,F(d_h\,y,s,d_h^{-1}p,d_h^{-2}M)
$$
A simple computation shows that the function $v(y)=u(d_h\,y)$, where $ u$ satisfies (\ref{subsol.f})  is  a viscosity solution of
$$
G(y,v(y),Dv(y),D^2v(y))\ge d_h^2\, f(d_hy)
$$
for $y \in d_h^{-1}\Omega \subset (0,1)\times \mathbb R^{n-1}$.

Since $G$ satisfies the same conditions as $F$ we may the apply the result for the case $d_h=1$ to $v(y)=u(x)$ and conclude the proof. \hfill$\Box$

\noindent We are now in position to state and give the simple proof based on the result of Theorem \ref{ABP-simple-thm} of ({\bf MP}) in narrow domains; that is, domains satisfying condition (\ref{slab-d}) with some small $d_h$, as stated in Theorem \ref{narrow-thm}.

\noindent {\bf Proof of Theorem \ref{narrow-thm}.} Let $u$ be as in $({\bf MP})$. As already observed $F(x,u,Du,D^2u) \ge 0$ implies $F(x,u^+,Du^+,D^2u^+) \ge 0$ and so, by $(SC)_{U,\Omega}$ ,
$$ F(x,0,Du^+,D^2u^+)\ge F(x,u^+,Du^+,D^2u^+)-cu^+\ge -cu^+.$$
By the assumptions on $u(x)$ and $\frac{c(x)}{\lambda(x)}$, bounded above, Theorem (\ref{ABP-simple-thm}) applies with $f=-c\,u^+$
yielding
$$\sup_{\overline\Omega} u \le \frac{e^{1+d\,\Gamma }}{1+d\,\Gamma} \,d^2\,\sup_{\Omega}\frac{c(x)u^+(x)}{\lambda(x)} \le  \frac{e^{1+d\,\Gamma }}{1+d\,\Gamma} \,d^2\,K  \sup_{\overline\Omega} u^+  $$
since $u^+= 0$ on $\partial \Omega$. From this estimate the statement immediately follows if $d$ is small enough.  \hfill$\Box$

\begin{oss}
{\rm It is worth observing that the same conclusion holds true for any $d$ such that $d_h \le d$ for all $h=1,\dots,k$, provided that  $\sup_\Omega c(x)$ is small enough. This fact is well-known for the case of uniformly elliptic linear operators of  the form ${\rm Tr }(A(x) D^2 u) +c(x) u$ if $c^+$ is small enough with respect to the ellipticity constant of the matrix $A$; see for instance \cite{Bu}, \cite{Vit}.}
\end{oss}

\section{Phragm\`en-Lindel\"of principles} \label{Phragmen-Lindelof} 

\noindent As in the proof of Theorem   \ref{ABP-simple-thm}, here we will use the fact that subsolutions of the equation $F(x,u,Du,D^2u) =0$ are in turn solutions of the differential inequality  
\begin{equation}\label{subsol.0}
F(x,0,Du^+,D^2u^+)\ge0
\end{equation}
{\bf Proof of Theorem \ref{Phragmen}.}  

We may suppose $\nu^1,\dots,\nu^k$ are the positive directions $E^1,\dots,E^k$ along $x_1,\dots,x_k$ so that the orthogonal subspace is generated by the positive directions along $x_{k+1},\dots,x_n$, and $P$, $Q$ are as in the proof of Theorem 1. \\
Suppose that $\Omega$ is contained in a $(n-k)$-infinite cylinder as defined in (\ref{slab-d}) with $\nu^h=E^h$, $h=1,\dots,k$. \\
Let us fix $h$. Our aim is to prove that for sufficiently small $d_h>0$ and viscosity solutions $u \in USC(\overline \Omega)$ of the differential inequality (\ref{subsol.0}) such that  $u(x) \le 0$ on $\partial\Omega$ and $u(x)=O(e^{\beta_0|x|})$ for a suitable $\beta_0>0$, we have in turn $u \le 0$ in $\Omega$.\\
Let $x=(y,z) \in \mathbb R^k\times \mathbb R^{n-k}$, $r=|z|=\left(x_{k+1}^2+\dots+x_n^2\right)^{\frac12}$ and set 
$$v(x)=\sin\alpha y_h\,e^{\beta \varphi(r)},\quad \varphi(r)=\sqrt{r^2+1}$$ 
where  $\beta>\beta_0$ and $\alpha > \beta$ will be chosen as follows\,. 
Computing $Dv(x)$ we get
$$
Dv(x)= \left(\begin{array}{cl}
&\, \vdots \\
& 0 \\
\alpha &\hskip-0.25cm\cos\alpha y_h \ e^{\beta \varphi(r)}\\
&0 \\
&\, \vdots  \vspace{0.2cm}\\
&\hskip-0.25cm v(x)\, e^{-\beta\varphi(r)} \frac{d}{dz}e^{\beta \varphi(r)}
\end{array}\right)
$$ 
where
$$e^{-\beta\varphi(r)}\left|\frac{d}{dr}e^{\beta \varphi(r)}\right|=\beta\varphi'(r)=\beta\,\frac{r}{\varphi(r)},$$
from which
\begin{equation}\label{Phragmen-eq1}
e^{-\beta\varphi}|Dv(x)|\le \alpha+\beta\,.
\end{equation}
Computing $D^2v(x^c)$, where $x^c=(\frac{\pi}{2\alpha},\dots,\frac{\pi}{2\alpha}, x'')$, we find
$$
 e^{-\beta \varphi(r)}D^2v(x^c)=\left(\begin{array}{cccc}
0 & 0& \dots&  0\\ 
\vdots &  -\alpha^2& \vdots & \vdots\\
0& \vdots&0 &   0\\ 
0 & 0 &\dots &  e^{-\beta \varphi(r)}D_{\hskip-0.05cm _{zz}}e^{\beta \varphi(r)}
\end{array}\right),$$
where
$e^{-\beta\varphi(r)}D_{\hskip-0.05cm _{zz}}e^{\beta \varphi(r)}$ is a positive definite $(n-k)\times (n-k)$ real matrix, having eigenvalues 
$$\beta\left(\beta (\varphi'(r))^2+ \varphi''(r) \right) = \beta\left(\beta\,\frac{ r^2}{\varphi^2(r)}+\frac1{\varphi^3(r)}\right) \;\; \mbox{\text and}\;\; \frac{\beta \varphi'(r)}{r}= \frac{\beta}{\varphi(r)}$$ 
with multiplicity $1$ and $n-k-1$, respectively (here $k<n$).
It follows that
\begin{equation}\label{Phragmen-eq2}
e^{-\beta\varphi(r)}D^2v(x^c) \le -\alpha^2 P_h + \beta(\beta+1)Q,
\end{equation}
where $P_h$ is the orthogonal projection on $x_h$. \\
For our purpose, we will choose $\alpha>0$ large enough in order that
\begin{equation}\label{Phragmen-eq3}
-\frac12\,\alpha^2 + 2\rho \beta(\beta+1) + \Gamma (\alpha +\beta) \le 0\,,
\end{equation}
where $\rho>0$ is an upper bound for $\Lambda(x)/\lambda(x)$, by assumptions.\\
Note that from (\ref{Phragmen-eq2}) by continuity there exists $d_0\in (0, \frac{\pi}{\alpha})$ such that for $x_h \in (\frac\pi{2\alpha}-\frac {d_0}2,\frac\pi{2\alpha}+\frac {d_0}2)$ we have 
\begin{equation}\label{Phragmen-eq4}
e^{-\beta\varphi(r)}D^2v(x) \le -\frac12\,\alpha^2 P_h + 2\beta(\beta+1)Q
\end{equation}
for all $x \in S\equiv  (0,d_1)\times \dots\times (\frac\pi{2\alpha}-\frac {d}2,\frac\pi{2\alpha}+\frac {d}2)\times \dots \times (0,d_k) \times \mathbb R^{n-k}$.

\noindent  Next we set 
$$w(x)=u^+(x)-c_Rv(x), \ \ x \in \Omega_R=\Omega \cap B_R(0)$$
where
$$c_R = \frac{\sup_{x \in \partial\Omega_R}u^+(x)} {e^{\beta R}\cos (\alpha \frac {d} 2)}.$$
Since $\Omega \subset S$, assuming $u \le 0$ on $\partial \Omega$, we have $w \le 0$ on $\partial \Omega_R$.\\

On the other hand, using  $(SC)_{U,\Omega}$ and (\ref{Phragmen-eq4}), from (\ref{Phragmen-eq3}) we get
\begin{equation*}\begin{split}
F(x,0,Dw,D^2w) &=\, F(x,0,Du^+-c_RDv,D^2u^+-c_RD^2v) \\
&\ge\, F(x,0,Du^+-c_RDv,D^2u^+ +c_R e^{\beta\varphi}(\tfrac12\,\alpha^2 P_h - 2\beta(\beta+1)Q)) \\
&\ge\, F(x,0,Du^+,D^2u^+) \\
&+\,c_Re^{\beta\varphi} \left(\tfrac12\,\lambda \alpha^2 -  2\Lambda\beta(\beta+1)-\gamma(\alpha+\beta)\right) \ge 0\,.
\end{split}
\end{equation*}
Hence Theorem  \ref{MP-slab-thm} yields $w \le 0$ in the bounded domain $\Omega_R$, namely
$$
u(x) \le c_Rv(x) = \frac{\sup_{\partial\Omega_R}u^+} {e^{\beta R}\cos (\alpha \frac {d} 2)} \, v(x)\,.
$$
Finally, consider an arbitrary $x \in \Omega$, and choose $R>0$ big enough in order that $x \in \Omega_R$. Letting $R \to \infty$ in the above, since $\sup_{\partial\Omega_R}u^+=O(e^{\beta_0R})$ and $\beta>\beta_0$, then $u(x) \le 0$.\\
The counterpart, that is the validity of the Phragm\`en-Lindel\"of principle with a suitable exponential growth for a fixed thickness $d_0>0$ of bounded orthogonal sections, is proved analogously.   \hfill$\Box$

\section*{Acknowledgements}
The authors would like to thank the referees for helpful comments and suggestions, which improved the presentation of the paper.\\
The authors are also grateful to Gruppo Nazionale per l'Analisi Matematica, la Probabilit\'a e le loro Applicazioni (GNAMPA) of Istituto Nazionale di Alta Matematica (INdAM).\\

\end{document}